\documentclass[11pt]{amsart}
\usepackage{amssymb}
\numberwithin{equation}{section}

\newtheorem{theorem}{Theorem}[section]
\newtheorem{Theorem}[theorem]{Theorem}
\newtheorem{corollary}[theorem]{Corollary}
\newtheorem{proposition}[theorem]{Proposition}

\theoremstyle{definition}
\newtheorem{definition}[theorem]{Definition}
\theoremstyle{remark}
\newtheorem{remark}[theorem]{Remark}
\newtheorem{example}[theorem]{Example}

\newcommand{\N}{\mathbb{N}}
\newcommand{\R}{\mathbb{R}}
\newcommand{\C}{\mathbb{C}}
\newcommand{\Z}{\mathbb{Z}}
\newcommand{\Q}{\mathbb{Q}}
\newcommand{\bbT}{\mathbb{T}}
\newcommand{\T}{\mathbb{T}}
\newcommand\lie[1]{\mathfrak{#1}}

\newcommand{\fg}{\lie{g}}

\def	\inv	{^{-1}}

\newcommand{\Mg}[1]{M_{[ #1, \infty)}} 
\newcommand{\Ml}[1]{M_{(-\infty, #1]}}

\begin{document}

\title{Contact cuts}

\author{Eugene Lerman}
\address{Department of
Mathematics, University of Illinois, Urbana, IL 61801}
\email{lerman@math.uiuc.edu}

\thanks{Supported by the NSF grant DMS-980305.}
\date{\today}

\begin{abstract}
We describe a contact analog of the symplectic cut construction
\cite{l.cut}.  As an application we show that the group of
contactomorphisms for certain overtwisted contact structures on
lens spaces  contains countably many non-conjugate two tori.
\end{abstract}

\maketitle

\section{Introduction}
In this paper we introduce the notion of cuts in the contact category
and use it to study the group of contactomorphisms of certain overtwisted contact structures on 3-manifolds.

Recall that a contact form on a $2n+1$ dimensional manifold $M$ is a
1-form $\alpha$ such that 
$$
\alpha \wedge (d\alpha )^n \not = 0.
$$ 
Thus a contact form defines an orientation on $M$. Note that if $n$ is
odd then $\alpha$ and $-\alpha$ define the same orientation. It is not
hard to see that if $f$ is a nowhere vanishing function on the
manifold $M$ and $\alpha $ is a contact form then $f\alpha$ is also a
contact form.  The kernels of $\alpha$ and $f\alpha$ are, of course,
the same.  A {\bf co-oriented contact structure} on a contact manifold $M$
is a subbundle $\xi$ of the tangent bundle $TM$ which is given as the
kernel of a contact 1-form.  A {\bf contact structure} $\xi$ on $M$ is
a subbundle of $TM$ which is locally, but not necessarily globally,
the kernel of a contact form.

Throughout the paper $\alpha$ will always denote a contact form and
$\xi$ will always denote a contact structure.  Whenever convenient we
will refer to a pair $(M, \alpha)$ or to a pair $(M, \xi)$ as a
contact manifold.

A diffeomorphism of a contact manifold $(M, \xi)$ preserving the
contact structure $\xi$ is called a {\bf contactomorphism}.

A contact 3-manifold  $(M, \xi)$   is {\bf symplecticly 
fillable} if there exists a compact symplectic manifold $(W, \omega)$
with boundary such that
\begin{enumerate}
\item $\partial W = M$,
\item $\omega|_\xi$ does not vanish, and
\item the orientation of $M$ defined by $\xi$ agrees with the orientation 
of $M$ as the boundary of the symplecticly oriented manifold $W$. 
\end{enumerate}
One often says in this case that $\xi$ is fillable.

A contact structure $\xi$ on a 3-manifold $M$ is {\bf overtwisted} if
there is an embedded 2-disk $D\subset M$ such that the boundary
$\partial D$ is tangent to $\xi$ and $D$ is transverse to $\xi$ along
$\partial D$.  A contact structure on a 3-manifold which is not
overtwisted is called {\bf tight}.  A theorem due to Gromov and
Eliashberg shoes that a fillable contact structure is tight
\cite[Theorem~3.2]{El3} (see also \cite{Gromov, El2}).  Thus the
standard contact structure on $S^3$ induced by its embedding in $\C^2$
as the unit sphere (the contact distribution consists of complex lines
tangent to the sphere) is tight.

Eliashberg showed \cite{Eliash} that on closed 3-manifolds overtwisted
contact structures are classified by the homotopy classes of the
corresponding plane fields.  We will exploit this property of
overtwisted contact structures to show that for a  class of such
structures the corresponding groups of contactomorphisms contain
countably many non-conjugate 2-tori.


\section{Cuts}  

\subsection{Topological cuts}

We start with the following observation.

\begin{proposition}\label{prop1}
Consider a smooth  action $(\lambda, m)\mapsto \lambda \cdot m$ of a circle
$S^1$ on a  smooth manifold $M$.  Let $f:M\to \R$ be an
$S^1$ invariant function.  Suppose that $a\in \R$ is a regular value
of $f$ and suppose that the action of $S^1$ on $f\inv (a)$ is free.
Then the topological space
$$
\Mg{a}:= \{ m \in M\mid f (m) \in [a, \infty) \}/\sim , 
$$
where, for $m\not = m'$, $m\sim m'$ if and only if 
\begin{enumerate}
\item $f (m) = f (m') = a$ and
\item $m = \lambda \cdot m'$ for some $\lambda \in S^1$, 
\end{enumerate}
is a smooth manifold.  The set $\{m\in M \mid f(m)> a\}$ is open and
dense in $\Mg{a}$ and the difference $\Mg{a}\smallsetminus \{m\in M
\mid f(m)> a\}$ is diffeomorphic to $f\inv (a)/S^1$.
\end{proposition}

We will refer to $\Mg{a}$ as {\bf the cut of $M$ with respect to the
ray $[a, \infty)$} (with the action and the function $f$ understood).

\begin{proof} 
The circle $S^1 \subset \C$ acts on the product manifold $M\times \C$
by $\lambda \cdot (m, z) = (\lambda \cdot m, \lambda\inv z)$.  The
function $\Psi (m, z) = f(m) - |z|^2$ is $S^1$ invariant.

The level set $\Psi\inv (a)$ is a (set theoretic) disjoint union of
two manifolds 
$$
\Psi\inv ( a)  =\left\{ (m,z): f(m) > a, \quad z = e^{i\theta}
\sqrt{f(m) - a} \right\} \bigsqcup \{(m,0): f (m)= a \}.
$$ 
The points of the first manifold are regular points of $\Psi$
because they are regular points of the function $(m, z) \mapsto
-|z|^2$.  The points of the second manifold are regular points of
$\Psi$ by assumption on $f$.  Moreover, by assumption on $f$ and by
construction the circle acts freely on $\Psi\inv ( a)$.  Therefore the
quotient $\Psi\inv ( a)/S^1$ is a smooth manifold.

The composition of $\sigma : \{ m\in M \mid f(m)\geq a\} \to \Psi \inv
(a)$, $\sigma (m) = (m, \sqrt{ f(m)- a})$ with the orbit map $\Psi
\inv (a) \to \Psi \inv (a)/S^1$ descends to a homeomorphism
$$
\bar{\sigma} : \Mg{a} \to \Psi \inv (a)/S^1.
$$ This gives the cut $\Mg{a}$ the structure of a smooth manifold.
Moreover, $\bar{\sigma}$ restricted to $\{ f(m) > a\}$ is an open
embedding and $\Psi\inv (a)/S^1 \smallsetminus \bar{\sigma} (\{ f(m) >
a\})$ is $f\inv (a)/S^1$.
\end{proof}

Completely analogously one defines 
$$
\Ml{a}:= \{ m \in M\mid f (m) \in  (-\infty, a] \}/\sim , 
$$
{\bf the cut of $M$ with respect to the ray $(-\infty, a]$}.  It is
also a manifold.

\begin{remark}
The action of $S^1$ on the manifold $M$ descends to an action of $S^1$
on the cut $\Mg{a}$.  The restriction $f|_{f\inv ([a, \infty))}$
descends to an $S^1$ invariant function $\bar{f}$ on the cut.

More generally if a group $G$ acts on $M$ preserving the function $f$
and commuting with the action of $S^1$, then the cut $\Mg{a}$ has a
naturally induced action of $G$ which again commutes with the action
of $S^1$ and preserves $\bar{f}$.
\end{remark}

It is not hard to see that one does not need a global circle action
for the construction in Proposition~\ref{prop1} to work.

\begin{proposition} \label{prop2}
Suppose $M$ is a manifold with boundary $\partial M$.  Suppose
further that $\partial M$ is a principal $S^1$ bundle.  Let $X = M/ \sim$,
where, for $m\not = m'$, $m\sim m'$ if and only if 
\begin{enumerate}
\item $m, m' \in \partial M$ and 
\item $m = \lambda \cdot m'$ for some $\lambda \in S^1$, 
\end{enumerate}
Then $X$ is a smooth manifold, $\partial M/S^1$ is a submanifold of
$X$, and $X \smallsetminus (\partial M/S^1)$ is diffeomorphic to $M
\smallsetminus \partial M$.
\end{proposition}

\begin{proof}
We may assume that a neighborhood  of the boundary $\partial M$ in $M$ is
diffeomorphic to $\partial M\times [0, 1)$.  Then $M = (M\smallsetminus
\partial M)\cup_{\partial M \times (0,1)} ( \partial M \times [0, 1))$.

The circle action on $\partial M$ extends trivially to a circle action
on $\partial M\times (-1, 1)$ making the projection map $f: \partial
M\times (-1, 1) \to (-1, 1) $ $S^1$-invariant.  By
Proposition~\ref{prop1}, $\left( \partial M \times (-1, 1)\right) _{[0,
\infty)}$ is a smooth manifold.  Moreover, it is easy to see that 
$$ X=(M\smallsetminus \partial M)\cup_{\partial M\times (0,1)}\left(
\partial M \times (-1, 1)\right) _{[0, \infty)}, 
$$ 
hence is smooth as well.
\end{proof}

\begin{example}
Let $M= \T^2 \times [0, 1]$ where $\T^2$ is the standard torus $\R^2
/2\pi \Z^2$.  Then $\partial M = \T^2 \times \{0, 1\}$.  Let $S^1$
act on $\partial M $ as follows: $S^1$ acts on $\T^2 \times \{0\}$ by $\lambda
\cdot (x, y) = (x + \lambda \mod 2\pi \Z, y)$ and  on $\T^2
\times \{1\}$ by $\lambda \cdot (x, y) = (x, y+ \lambda \mod 2\pi \Z)$.
Then $X$ is $S^3$.
\end{example}

\subsection{Symplectic cuts}

Next suppose that additionally the manifold $M$ possesses an $S^1$
invariant symplectic form $\omega$ such that the function $f$ is a
moment map for the action of $S^1$.  Then the cut $\Mg{a}$ is a
symplectic manifold.  More precisely:

\begin{Theorem}\label{symp_cut}
Let $(M, \omega)$ be a symplectic manifold with a Hamiltonian action
of a circle $S^1$, let $f:M \to \R$ denote a corresponding moment map.
Suppose that $S^1$ acts freely on the level set $f\inv (a)$ for some
$a\in \R$ (so that $a$ is a regular value of $f$).  Then the cut of
$M$ with respect to the ray $[a, \infty)$ is naturally a symplectic
manifold.  Moreover, the natural embedding of the reduced space $M_a
:= f\inv (a)/S^1$ into $\Mg{a}$ is symplectic, and the complement
$\Mg{a} \smallsetminus M_a$ is symplectomorphic to the open subset
$\{m\in M \mid f(m)> a\}$ of $(M, \omega)$.
\end{Theorem}

\begin{proof}  (cf. \cite{l.cut})
Consider the symplectic product $(M\times \C, \omega + \sqrt{-1}dz
\wedge d \bar{z})$.  The map $\Psi (m, z) = f(m) - |z|^2$ is a moment
map for an action of $S^1$ on $M\times \C$.

Arguing as in Proposition~\ref{prop1} we see that $a$ is a regular
value of $\Psi$ and that the reduced space $\Psi\inv ( a) /S^1$ is the
cut $\Mg{a}$.

 The pullback of the symplectic form $\omega + \sqrt{-1}dz
\wedge d \bar{z}$  by the embedding
$$
\sigma : \{ m\in M : f(m)> a\} \hookrightarrow \Psi\inv (a), \qquad
\sigma (m) = (m, \sqrt{f(m) -a })
$$
is $\omega$.  Consequently the induced embedding $\bar{\sigma}: \{ m\in M :
f(m)> a\} \hookrightarrow \Psi\inv (a)/S^1 = \Mg{a}$ is symplectic.
Similarly one checks that the natural embedding $M_a\hookrightarrow
\Mg{a}$ is symplectic  as well.

\end{proof}
\begin{remark}
More generally the construction can be carried out for Hamiltonian
torus actions.  A ray would then be replaced by a simple rational
polyhedral cone in the dual of the Lie algebra of the torus with the
``Delzant condition'' on the edges; see \cite{l.cut, toplogy}.
\end{remark}

An analog of Proposition~\ref{prop2} holds as well:

\begin{proposition}
Suppose $(M, \omega)$ is a symplectic manifold with boundary $P
=\partial M$.  Suppose further that $P$ is a principal $S^1$ bundle,
that $\omega|_P$ is $S^1$ invariant and that the kernel of $\omega|_P$
is precisely the vertical bundle of $P\to P/S^1$.  Let $X = M/ \sim$,
where, for $m\not = m'$, $m\sim m'$ if and only if
\begin{enumerate}
\item $m, m' \in P$ and 
\item $m = \lambda \cdot m'$ for some $\lambda \in S^1$, 
\end{enumerate}
Then $X$ is a symplectic manifold, $P/S^1$ is a symplectic
submanifold of $X$, and $X \smallsetminus (P/S^1)$ is symplectomorphic to
$M\smallsetminus P$.
\end{proposition}

\begin{proof}
By Proposition~\ref{prop2} $X$ is a smooth manifold, $P/S^1$ is a
submanifold and $X \smallsetminus (P/S^1)$ is diffeomorphic 
to $M\smallsetminus P$.  

We now assume for simplicity that $P$ is connected.  Otherwise we can argue
connected component by connected component.

By the equivariant coisotropic embedding theorem the product $P\times
\R$ carries and $S^1$-invariant closed 2-form $\tilde{\omega}$ ($S^1$
acts on $P\times \R $ by $\lambda \cdot (p, t) = (\lambda \cdot p, t)$)
such that
\begin{enumerate}
\item $\tilde {\omega}|_{P\times \{0\}} = \omega |_P$
\item the $S^1$ action on $(P\times \R, \tilde {\omega}) $ is 
Hamiltonian with  a moment map $f(p,t) = t$.
\end{enumerate}
Moreover there is an open $S^1$ equivariant embedding $\psi$ of a
neighborhood $U$ of $P$ in $M$ into $P\times \R$ such that
\begin{enumerate}
\item $\psi (p) = (p, 0)$ for all $p\in P$,
\item $f\circ \psi (m) \geq 0$ for all $m\in U$ and 
\item $\psi ^* \tilde {\omega} = \omega$.
\end{enumerate}
In particular $\tilde{\omega}$ is non-degenerate near $P\times\{0\}$.

We have $M = (M\smallsetminus P) \cup_{U\cap (M\smallsetminus P) }U$.
By Theorem~\ref{symp_cut} $U_{[0, \infty)}$ is a symplectic manifold,
the embedding of $P/S^1$ into $U_{[0, \infty)}$ is symplectic and the
difference $U_{[0, \infty)} \smallsetminus (P/S^1)$ is
symplectomorphic to $U\smallsetminus P$.  Therefore $X=
(M\smallsetminus P) \cup_{(U\smallsetminus P) }U_{[0, \infty)}$ is a
symplectic manifold with the desired properties.
\end{proof}


\subsection{Contact cuts} 
We start by digressing on group actions on contact manifolds and
follow the digression by recalling the definition of the moment map
for a group action preserving a contact form \cite{?, Geiges}.

\begin{proposition}
Let $M$ be a paracompact manifold with a contact form $\alpha$, let
$\xi = \ker \alpha$ be the corresponding contact structure.
Suppose a Lie group $G$ acts properly on $M$, that is, suppose the map
$G\times M \to M\times M$, $(g, m) \mapsto (g\cdot m, m)$ is proper.
Suppose further that the action of $G$ preserves the contact structure
$\xi$.  

Then there exists a $G$-invariant contact form $\bar{\alpha}$ with
$\ker \bar{\alpha} =\xi$.
\end{proposition}

\begin{proof}
The argument is an adaptation of Palais's proof of the existence of
invariant Riemannian metrics on manifolds with proper group actions
\cite{Palais}.

Suppose first that the group $G$ is compact.  Then there is on $G$ a
bi-invariant measure $dg$ normalized so that $\int _G dg = 1$.
We then define $\bar{\alpha}$ to be the average of $\alpha$:
$$
\bar{\alpha}_x := \int _G (g^* \alpha)_x\, dg
$$
for all $x\in M$.

Now we drop the compactness assumption.  None the less, for every
point $x\in M$ there exists a slice $S$ to the action of $G$.  That
is, $S$ is $H= \text{stab} (x)$ -invariant embedded submanifold of $M$
such that the union $G\cdot S$ of $G$-orbits through the points of $S$
is open and such that for every $s\in S$ the orbit $G\cdot s$
intersects $S$ in a single $H$-orbit (see \cite{Palais}).  Note that
the stabilizer $H$ of $x$ is compact because the action is proper.

Since the contact structure $\xi$ is $G$-invariant its annihilator
$\xi^\circ \subset T^*M $ is a $G$-invariant line subbundle of the
cotangent bundle.  Hence the restriction $\xi^\circ |_{G\cdot S}$ is
completely determined by the restriction $\xi^\circ |_{S}$.  Moreover, any
$G$-invariant section $\sigma : G\cdot S \to \xi^\circ|_{G\cdot S}$ is
completely determined by its values on $S$: 
$$
\sigma (g\cdot s) = g\cdot \sigma (s)
$$
for all $g\in G$, $s\in S$.  

A contact form $\alpha$ is a nowhere vanishing section of $\xi^\circ$.
Given a slice $S$ we produce an $H$-invariant section $\bar{\alpha}$
of $\xi^\circ|_S$ by averaging $\alpha|_S$ over $H$.  We then extend
$\bar{\alpha}$ to the open set $G\cdot S$ by the formula
$\bar{\alpha}_{g\cdot s} = g\cdot \bar{\alpha}_s$.  The form
$\bar{\alpha}$ is well-defined on $G\cdot S$ because $(G\cdot s) \cap S
= H \cdot s$ for any $s\in S$.

Next cover $M$ by open sets $U_\beta$ of the form $G\cdot S_\beta$
where $S_\beta $ are slices.  In the proof of
\cite[Theorem~4.3.1]{Palais} Palais showed that this cover may be
chosen to be locally finite and that there exist a $G$-invariant
partition of unity $\{\rho _\beta\}$ subordinate to the cover.  The form
$$
\bar{\alpha} := \sum \rho_\beta \, \overline{ \alpha |_{G\cdot S_\beta}}
$$ 
is the desired invariant contact form.  Here $\overline{ \alpha
|_{G\cdot S_\beta}}$ is the average of $\alpha |_{G\cdot S_\beta}$.
\end{proof}

\begin{definition}
Let $M$ be a manifold with a contact one-form $\alpha$.  Suppose a Lie
group $G$ acts properly on $M$ and preserves the contact distribution
$\ker \alpha$.  By averaging if necessary (see above) we may assume
that $\alpha$ is $G$-invariant.  We define the corresponding {\bf
moment map }$\Phi :M\to
\fg^*$ ($\fg^*$ is the dual of the Lie algebra $\fg$ of $G$) by the equation
$$
\langle \Phi, \eta \rangle = \alpha (\eta_M)
$$ 
for all $\eta\in \fg$.  Here $\langle \cdot,
\cdot \rangle : \fg^* \times \fg \to \R$ is the standard pairing, and
$\eta_M$ denotes the vector field on $M$ induced by $\eta \in \fg$.
\end{definition}

The following result is a slight generalization of Theorem~6 in \cite{Geiges}. 
\begin{theorem}
Suppose $(M, \alpha)$ is a contact manifold with a proper action of a
Lie group $G$ preserving the contact form $\alpha$.  Suppose 0 is a
regular value of the corresponding moment map $\Phi: M\to \fg^*$.
Then $\alpha |_{\Phi\inv (0)}$ descends to a contact form $\alpha _0$
on the orbifold $M_0 := \Phi\inv (0)/G$.
\end{theorem}

\begin{proof}
The proof is identical to the proof of Theorem~6 in \cite{Geiges}.
The main idea of the proof is that a point $x$ lies in the zero level
set of the moment map if and only if the orbit $G\cdot x$ is tangent
to $\ker \alpha$; hence $\alpha |_{\Phi\inv (0)}$ descends to a 1-form
$\alpha_0$ on $M_0$.
\end{proof}
We will refer to the pair $(M_0, \alpha_0)$ as the contact quotient
of $(M, \alpha)$ or as the reduced space.

\begin{theorem}\label{contact_cut}
Let $(M, \alpha)$ be a contact manifold with an action of $S^1$
preserving $\alpha$ and let $ f$ denote the corresponding the moment
map.  Suppose that $S^1$ acts freely on the zero level set $f\inv
(0)$.  Then the cut $\Mg{0}$ of $M$ is naturally a contact manifold.
Moreover, the natural embedding of the reduced space $M_0 := f\inv
(0)/S^1$ into $\Mg{0}$ is contact and the complement $\Mg{0}
\smallsetminus M_0$ is contactomorphic to the open subset $\{m\in M
\mid f(m)> 0\}$ of $(M,\alpha)$.
\end{theorem}
\begin{proof}
Consider the contact manifold $(M\times \C, \alpha + \frac{1}{2}
\sqrt{-1} (zd\bar{z} - \bar{z} dz))$ with the circle action $\lambda
\cdot (m,z) = (\lambda \cdot m, \lambda \inv z)$.  The map 
$\Psi(m, z) = f(m) - |z|^2$ is the corresponding moment map. 
Arguing as in Proposition~\ref{prop1} we see that $0$ is a regular
value of $\Psi$ and that the reduced space $\Psi\inv (0) /S^1$ is the
cut $\Mg{0}$.

The pullback of the contact form $\alpha +
\frac{1}{2}\sqrt{-1} (zd\bar{z} - \bar{z} dz)$ by the embedding
$$
\sigma : \{ m\in M : f(m)> 0\} \hookrightarrow \Psi\inv (0), \qquad
\sigma (m) = (m, \sqrt{f(m)})
$$ is $\alpha$.  Consequently the induced embedding $\bar{\sigma}: \{
m\in M : f(m)> 0\} \hookrightarrow \Psi\inv (0)/S^1 = \Mg{0}$ is
contact.  Similarly one checks that the natural embedding
$M_0\hookrightarrow
\Mg{0}$ is contact as well.
\end{proof}

\begin{example}\label{ex1}
Consider the three-manifold $M = S^1\times S^1\times\R$ with
coordinates $\theta_1$, $\theta _2$ and $t$.  The form $\alpha = \cos
t \,d\theta _1 + \sin t \,d\theta _2$ on $M$ is contact.  The two
torus $\bbT^2 = S^1\times S^1$ acts on $M$ by 
$$ 
(a, b)\cdot
(\theta_1, \theta _2, t) = (\theta_1 +a , \theta _2 + b, t) 
$$ 
for all $(a, b)\in S^1\times S^1$, $(\theta_1, \theta _2, t) \in
S^1\times S^1\times\R$ (here we think of $S^1$ as $\R/2\pi \Z$).  The
action preserves the contact form $\alpha$, and the corresponding
moment map $f= (f_1, f_2) : M \to \R^2$ is given by 
$$ 
f(\theta_1, \theta _2, t) = (\cos t, \sin t).  
$$ 
Let's cut $M$ with respect to $[0, \infty)$ using the second component
$f_2 = \sin t$ of $f$.  Since 
$$ 
f_2 \inv ([0, \infty)) =
\{(\theta_1, \theta _2, t) : t\in [2\pi n, \pi (2n+1)], n \in \Z\},
$$
the cut is $\bigsqcup_{n \in Z} \{(\theta_1,
\theta _2, t) : t\in [2\pi n, \pi (2n+1)]\}/\sim$, which is the
disjoint union of countably many copies of $S^1\times S^2$.  We will
see later that the contact structure that we have constructed on
$S^1\times S^2$ is fillable.
\end{example}

\begin{remark}
As in the symplectic case, if the contact manifold carries a torus
action then we can define cuts with respect to a simple rational
polyhedral cone in the dual of the Lie algebra of the torus.  One can
prove directly that if the moment map is transverse to all faces of
the cone then the resulting cut space is a contact orbifold.
Alternatively, one can apply Theorem~\ref{contact_cut} above and
contact reduction in stages.
\end{remark}

\begin{remark}
The construction of a contact cut does not depend on a choice of a
contact form.  Suppose we change the contact form by multiplying it by
a positive invariant function.  The corresponding moment map gets
multiplied by the same function.  Consequently the contact form on the
cut space gets multiplied by a positive function.

Alternatively, since contact reduction can be defined without any
reference to moment maps \cite{Severa}, contact cuts may also be
defined without moment maps: we define it to be the contact reduction
of $M\times \C$.
\end{remark}

It is important to note that we have a contact analog of
Proposition~\ref{prop2}.

\begin{proposition}\label{cont_prop2}
Suppose $(\tilde M, \alpha)$ is a contact manifold, $M$ is a manifold
with boundary of the same dimension as $\tilde{M}$ embedded in
$\tilde{M}$.  Suppose further that there is a neighborhood $U$ in
$\tilde M$ of the boundary $\partial M$ and  
a free $S^1$ action on $U$ preserving $\alpha$ such
that the corresponding moment map $f:U \to \R$ satisfies
\begin{enumerate}
\item $f\inv (0) = \partial M $ and
\item $f\inv ([0, \infty)) = U \cap M$.
\end{enumerate}
 
Let $X = M/ \sim$, where, for $m\not = m'$, $m\sim m'$ if and only if
\begin{enumerate}
\item $m, m' \in \partial M$ and 
\item $m = \lambda \cdot m'$ for some $\lambda \in S^1$, 
\end{enumerate}
Then $X$ is a contact manifold, $\partial M/S^1$ is a contact
submanifold of $X$, and $X \smallsetminus (\partial M/S^1)$ is
contactomorphic to $M\smallsetminus \partial M$.
\end{proposition}

\begin{proof}
By Proposition~\ref{prop2} $X$ is a smooth manifold and $X
\smallsetminus (\partial M/S^1)$ is diffeomorphic to $M\smallsetminus
\partial M$.  We would like to show that $X$ is contact.

We have $M = (M\smallsetminus \partial M)
\cup_{U\cap (M\smallsetminus \partial M) }
U \cap M $.  By Theorem~\ref{contact_cut} $U_{[0, \infty)}$ is a
contact manifold, the embedding of $\partial M/S^1 $ into $U_{[0,
\infty)}$ is contact and the difference $U_{[0, \infty)}
\smallsetminus (\partial M/S^1)$ is contactomorphic to $U\cap
(M\smallsetminus \partial M)$.  Therefore $X = (M\smallsetminus
\partial M) \cup_{U\cap (M\smallsetminus \partial M)}
U_{[0, \infty)}  $ is a contact manifold with the desired properties.
\end{proof}

\begin{example}\label{exS3}
Let $\tilde {M}$ denote the manifold $S^1\times S^1\times\R$ with
coordinates $\theta_1$, $\theta _2$ and $t$ and a contact form $\alpha
= \cos \frac{\pi}{2} t \,d\theta _1 + \sin \frac{\pi}{2}t \,d\theta
_2$.  The manifold $M=S^1\times S^1\times [0, 1]$ embeds
naturally into $\tilde M$.  The boundary of $M$ is $\partial M = S^1
\times S^1
\times \{0, 1\}$.  Let 
$U_1 = S^1\times S^1 \times (-\frac{1}{2}, \frac{1}{2})$, 
$U_2 = S^1\times S^1 \times (\frac{1}{2}, \frac{3}{2})$.  

Consider the $S^1$ action on $U= U_1 \cup U_2$ induced on $U_1$ by the
vector field $\frac{\partial }{\partial \theta_2}$ and on $U_2$ by
$\frac{\partial }{\partial \theta_1}$.  The corresponding moment map
is $f$ is given by $f|_{U_1} = \iota (\frac{\partial }{\partial
\theta_2}) \alpha = \sin\frac{\pi}{2} t$ and 
$f|_{U_2} = \iota (\frac{\partial }{\partial
\theta_1}) \alpha = \cos\frac{\pi}{2} t$.  The set $\{f \geq 0\}$ is 
$\{ (\theta_1, \theta_2, t)\in \tilde{M} \mid 0\leq t <
\frac{1}{2}\, \text{or} \,\frac{1}{2} < t\leq 1\}$.
The cut manifold $X$ is easily seen to be the 3-sphere $S^3$.

We will see below (example~\ref{exT^3}) that the induced contact
structure on $X$ is the standard contact structure on $S^3$.
\end{example}

\subsection{Cuts of hypersurfaces of contact type}
 
It will be useful for us to have a way of establishing that the
contact structure on a contact cut is fillable.  Recall that a
hypersurface $\Sigma$ in a symplectic manifold $(M, \omega)$ is of
{\bf contact type} if there exists on a neighborhood of $\Sigma$ a
vector field $X$ such that $X$ is transverse to $\Sigma$ and such that
the Lie derivative of $\omega$ with respect to $X$ is $\omega$: $L_X
\omega = \omega$.  It is well known that in this case the restriction
of the contraction $\iota (X)\omega$ to $\Sigma$ is a contact form.

Our first step is to note that the contact quotient of an invariant
hypersurface of contact type is again a hypersurface of contact type.
More precisely:

\begin{proposition}\label{prop3}
Suppose $(M, \omega)$ is a symplectic manifold with a Hamiltonian
action of $S^1$ and a corresponding moment map $f: M \to \R$.  Suppose
that $\Sigma$ is a hypersurface in $M$ which is preserved by the
action of $S^1$.  Moreover, assume that there is an $S^1$ invariant
vector field $X$ defined in an invariant neighborhood of $\Sigma$ such
that $X$ is transverse to $\Sigma$ and such that the Lie derivative of
$\omega$ with respect to $X$ is $\omega$.  Then the contact form
$\alpha := (\iota (X) \omega)|_\Sigma$ is $S^1$ invariant, and the
moment map $f_\Sigma$ for the action of $S^1$ on $(\Sigma, \alpha)$
is, up to an additive constant, the restriction of $f$ to $\Sigma$.
We may assume that the constant is zero.

Suppose further that $S^1$ acts freely on $f\inv (0)$.  Then $0$ is a
regular value of $f$ and of $f|_\Sigma$, and the contact quotient
$\Sigma_0$ is a hypersurface of contact type in the symplectic
quotient $M_0 := f\inv (0)/S^1$.
\end{proposition}

\begin{proof} 
The vector field $X$ descends to a vector field $X_0$ on a
neighborhood of $\Sigma_0$.  It is not hard to see that the reduced
symplectic form $\omega_0$ satisfies $L_{X_0} \omega _0 = \omega _0$
and that $(\iota (X_0) \omega_0)|_{\Sigma_0}$ is the reduced contact
form $\alpha_0$.
\end{proof}

\begin{corollary} \label{cor_cont_type}
Suppose $(M, \omega)$, $\Sigma \subset M$ and $f:M \to \R$ are as in
the proposition above. 

The cut $\Sigma _{[0, \infty)}$ is a hypersurface of contact type in
the cut $\Mg{0}$, and the cut contact form on $\Sigma _{[0, \infty)}$
is $(\iota (\bar{X})\omega _{[0,\infty)})|_{\Sigma _{[0, \infty)}}$,
where $\omega _{[0, \infty)}$ is the induced symplectic form on
$\Mg{0}$ and $\bar{X}$ is the vector field on a neighborhood of
$\Sigma _{[0, \infty)}$ induced by $X$.
\end{corollary}

\begin{proof}
The hypersurface $\Sigma \times \C$ in $M\times \C$ is of contact type.
The cut $\Sigma _{[0, \infty)}$ is the reduction of $\Sigma \times
\C$, and the cut $\Mg{0}$ is the reduction at zero of $M\times \C$.  Now apply 
Proposition~\ref{prop3}.
\end{proof}
\begin{example}\label{exT^3}
Consider the cotangent bundle $M= T^*\bbT^2$ of the standard two torus
with the standard symplectic form.  Denote the coordinates on $\bbT^2$
by $\theta_1$, $\theta_2$ and the corresponding coordinates on $T^*
\bbT^2$ by $\theta_1$, $\theta_2$, $p_1$, $p_2$.  In these coordinates
the cosphere bundle with respect to the flat metric is
$$
\Sigma = \{( \theta_1, \theta_2,p_1, p_2) : p_1^2 + p_2^2 = 1\}
\simeq \bbT^3.
$$
The contact form on $\Sigma$ is 
$$
\alpha = \sum p_i \,d\theta _i |_\Sigma.
$$
Consider the action of $S^1$ given by 
$$
a \cdot (\theta_1, \theta _2, p_1, p_2 ) = 
(\theta_1 , \theta _2 + a, p_1, p_2 )
$$ 
The map 
$$
f(\theta_1, \theta _2, p_1, p_2) =  p_2
$$
is a corresponding moment map.  It is not hard to see that $\Mg{0} =
\{( \theta_1, \theta_2,p_1, p_2) : p_2\geq 0\}/\sim $ is 
symplectomorphic to $ T^*S^1\times \C$, and that $\Sigma_{[0,\infty)}
\simeq S^1 \times S^2$.  Consequently the contact structure on
$S^1\times S^2$ that we obtained by cutting is fillable, hence tight
\cite{El4}  (cf. example~\ref{ex1}).

Next consider the action of $S^1$ on $T^*\T^2$ given by 
$$ 
b\cdot(\theta_1, \theta _2, p_1, p_2 ) = (\theta_1 +b , \theta _2 , p_1, p_2) 
$$ 
The map 
$$
 h(\theta_1, \theta _2, p_1, p_2) = p_1 
$$ 
is a corresponding moment map.  This action of $S^1$ commutes with the
first action of $S^1$.  Consequently it descends to a Hamiltonian
action on $\Mg{0} \simeq T^*S^1\times \C$.  The corresponding moment
map $\bar{h}: T^*S^1\times \C \to \R$ is given by $\bar{h} (\theta, p,
z) = p + |z|^2$ for all $(\theta, p, z) \in T^*S^1\times \C$.  If we
now cut again we obtain $ (T^*S^1\times \C)_{[0,\infty)} \simeq \C^2$.
The cut of the hypersurface is the standard $S^3$ in $\C^2$.  It
follows that the contact structure constructed in Example~\ref{exS3}
is the standard tight contact structure on the 3-sphere.
\end{example}

\subsection{Cuts and symplectization}

A contact manifold $(M, \alpha)$ embeds naturally as a hypersurface of
contact type into its {\bf symplectization} $(M\times \R, d (e^t
\alpha))$, where $t$ is the coordinate on $\R$.  Moreover, if  
there is an action of $S^1$ on $M$ preserving $\alpha$, then the
trivial extension of the action to $M\times \R$ preserves the
symplectic form $d(e^t \alpha)$.  The two moment maps $\Psi :M\times
\R
\to \R$ and $\Phi : M\to \R$ are related by the formula 
$$
\Psi (m, t) = - e^t \Phi (m)
$$ 
for all $(m,t) \in M\times \R$. Using this formula it is easy to see
that symplectization and reduction commute: 
$$ 
((M\times \R)_0, (d(e^t \alpha))_0) = (M_0 \times \R, d(e^t \alpha_0)) 
$$ 
where $(d(e^t \alpha))_0$ denotes the reduced symplectic form and
$\alpha_0$ denotes the reduced contact form.

Since cutting amounts to reduction, we see that cuts and symplectization
commute:

\begin{proposition}
Let $(M, \alpha)$ be a contact manifold with an action of $S^1$
preserving the contact form.  The the symplectization $\Mg{0} \times
\R$ of the cut $\Mg{0}$ is symplectomorphic to the cut of the
symplectization $(M\times \R)_{[0, \infty)}$.
\end{proposition}

\subsection{Contact gluing}

Recall that given a symplectic manifold $(M, \omega)$ with a function
$f: M \to \R$ generating a Hamiltonian circle action we obtain two cut
manifolds $\Mg{a}$ and $\Ml{a}$ for every regular value $a$ of $f$.
Conversely, the manifold $(M, \omega)$ can be reconstructed from the
cuts $\Mg{a}$ and $\Ml{a}$ by gluing them symplecticly along the
codimension two submanifold $M_a$ (see \cite{Gompf} for a precise
description of symplectic gluing).

Symplectic gluing has a counterpart in the contact category.  If $N$
is a contact submanifold of a contact manifold $(M, \alpha)$ then the
normal bundle $\nu$ of $N$ is symplectic, hence has well-defined Chern
classes.  Suppose now that a manifold $N$ is embedded as a codimension
two contact submanifold of two contact manifolds $(M_1, \alpha_1)$ and
$(M_2, \alpha_2)$ such that the first Chern classes of the
corresponding normal bundles $\nu_1$, $\nu_2$ are the negatives of
each other: $c_1 (\nu_1) = - c_1 (\nu_2)$.  Then one can glue $M_1$
and $M_2$ along $N$ to obtain a new contact manifold $M_1 \#_N
M_2$ with the following properties:
\begin{enumerate}
\item It contains a hypersurface $P$ diffeomorphic to the sphere bundle of 
$\nu_1$ (and of $\nu_2$).

\item $M_1 \#_N M_2 \smallsetminus P$ is contactomorphic to the disjoint  
union of $M_1 \smallsetminus N$ and $M_2 \smallsetminus N$. 
\end{enumerate}
This construction, with the additional assumption that $c_1 (\nu_i) =
0, \, i=1,2$, is described in detail in \cite{Geiges}.  The
construction also works in full generality.

Now if $(M, \alpha)$ is a contact manifold with a circle action
preserving the contact form $\alpha$ and the corresponding moment map
$f: M\to \R$, the contact cuts $\Mg{0}$ and $\Ml{0}$ can be glued
contactly along $M_0$ to reconstruct $M$ up to a contactomorphism.

\section{An application}

Recall that according to Eliashberg's classification of co-oriented
contact structures on on the 3-sphere there exists exactly one
overtwisted co-oriented contact structure $\zeta_0$ on $S^3$ which is
homotopic to the standard contact structure as a 2-plane field (cf.\
\cite{ABKLR}) .  Let us refer to $\zeta _0$ as the standard overtwisted
structure on $S^3$.

In this section we use contact cuts to produce countably many
non-conjugate 2-tori in the group of contactomorphisms of the standard
overtwisted contact structure on the 3-sphere.  We then remark that
for any rational number $p/q$ the same construction produces an
overtwisted contact structure $\zeta$ on a lens space $L_{p/q}$ with
the same property: the group of contactomorphisms of $(L_{p/q},
\zeta)$ contains countably many non-conjugate 2-tori.

\begin{theorem}\label{main}
Consider the standard overtwisted contact structure $\zeta_0$ on the
3-sphere $S^3$. The group of contactomorphisms of $(S^3, \zeta_0)$
contains countably many non-conjugate two-tori.  
\end{theorem}

\begin{proof}
As in example~\ref{exS3} let $\tilde {M}$ denote the manifold
$S^1\times S^1\times\R$ with coordinates $\theta_1$, $\theta _2$ and
$t$ and a contact form $\alpha = \cos t \,d\theta _1 + \sin t
\,d\theta _2$.  For each non-negative integer $k$ consider an embedding 
$\iota_k$ of $M = S^1\times S^1 \times [0, 1]$ into $\tilde{M}$:
$$
\iota_k (\theta_1, \theta_2, t) =  (\theta_1, \theta_2, (2\pi k + \pi/2)t).
$$ 
Let $\alpha _k = \iota_k ^* \alpha = \cos [(2 k + 1/2)\pi t] \,
d\theta _1 + \sin [(2 k + 1/2)\pi t] \, d\theta _2$.  Consider the
action of $S^1$ on the boundary  $\partial M$ generated by 
$\frac{\partial}{\partial \theta
_2}$ on $S^1\times S^1 \times \{0\}$ and by $\frac{\partial}{\partial
\theta _1}$ on $S^1\times S^1 \times \{1\}$.  Then as in example~\ref{exS3},
the cut manifold $X = M/\sim$ is $S^3$.  By
Proposition~\ref{cont_prop2} each of the contact forms $\alpha_k$
induces a contact from $\bar{\alpha}_k$ on $X$.

Note that by examples~\ref{exS3} and \ref{exT^3} the contact from
$\alpha_0$ defines the standard tight contact structure on $S^3$.
We now argue that
\begin{enumerate}
\item  All contact forms $\alpha _k$ are homotopic as non-vanishing one forms,
that is, for all $k, l\in \N$ there is a family of 1-forms
$\{\alpha^{k,l}_t\}_{0\leq t \leq 1} $ with $\alpha^{k,l}_0 =
\bar{\alpha}_k$, $\alpha^{k,l}_1 = \bar{\alpha}_l$, and
$(\alpha^{k,l}_t)_x \not = 0$ for all $x\in X$.  Hence the
corresponding contact structures are all homotopic as 2-plane fields.

\item  For $k>0$ the contact structures $\ker \bar{\alpha}_k$ are overtwisted.
Since overtwisted contact structures on 3-manifolds are classified by
the homotopy type of the corresponding 2-plane fields \cite{Eliash},
it would then follow that all the contact structures $\bar{\alpha}_k$
are equivalent provided $k>0$.

\item For any two  distinct integers $k, l \geq 1$,  there is no 
contactomorphism $\psi: (X, \bar{\alpha}_k) \to (X, \bar{\alpha}_l)$
which is $S^1\times S^1$ equivariant.
\end{enumerate}
Items (2) and (3) together  then prove the theorem.

To construct a path between $\bar{\alpha}_k$ and $\bar{\alpha}_l$
consisting of nowhere vanishing 1-forms, consider first the straight
line homotopy between $\alpha _k$ and $\alpha _l$ on $M$.  It fixes
the end points and so descends to a homotopy between $\bar{\alpha}_k$
and $\bar{\alpha}_l$ on $X$.  If the homotopy vanishes at a point
$x_0$ at a time $t_0$ correct it by pushing the path out in the $dt$
direction.

To see that the contact structures defined by $\bar{\alpha}_k$ ($k>0$)
are overtwisted, fix $c\in S^1$ and consider the subset $D^k_c \subset
X$ given by 
$$ 
D^k_c = \{ [\theta_1, \theta_2, t] \in X \mid \theta _1
= c, \quad 0\leq t \leq \frac{2\pi}{4k+1} \},
$$
where $[\theta_1, \theta_2, t]$ is the class of $(\theta_1, \theta_2,
t) \in M$.  The set $D^k_c$ is an overtwisted disk in 
$(X, \bar{\alpha}_k)$.

To prove (3) consider two contact connected manifolds $(N_1,
\alpha_1)$, $(N_2,\alpha_2)$ with an action of a Lie group $G$ preserving 
the contact forms.  Let $f_i: N_i \to \fg^*$ ($i=1,2$) denote the
corresponding moment maps.  If $\psi: (N_1, \alpha_1) \to (N_2,
\alpha_2)$ is a $G$-equivariant contactomorphism, then $\psi^*
\alpha_2 = h \alpha_1$ for some nowhere vanishing $G$-invariant
function $h$ on $N_1$.  Since $N_1$ is connected, $h$ is either always
positive or always negative.  Say $h>0$.  Then for any nonzero vector
$\xi \in \fg^*$ the preimages $f_i\inv (\{e^s \xi \mid s\in \R\})$ of
the ray through $\xi$ have the same number of connected components (if
$h<0$ then $\pi_0 (f_1\inv (\{e^s \xi \mid s\in \R\})) = \pi_0
(f_2\inv (\{-e^s \xi \mid s\in \R\}))$).

Now for any $k>0$ the map moment map $f_k: (X_k, \bar{\alpha}_k) \to
\R^2$ for the action of $S^1\times S^1$ is given by 
$$
f_k ([\theta_1, \theta_2, t])=
 (\cos ( \frac{4k+1}{2}\pi t), \sin( \frac{4k+1}{2}\pi t))
$$
Hence
$$
\pi_0 (f_k\inv (\{e^s (-1, 1) \mid s\in \R\}) =
\pi_0 (f_k\inv (\{e^s (1, -1) \mid s\in \R\}) = k.
$$
Therefore (3) follows and so does the theorem.
\end{proof}

Let us now turn our attention to lens spaces.  Recall that
the lens space $L_{p/q}$ is the manifold obtained by attaching two
solid tori $S^1 \times D^2$ together by a diffeomorphism sending a
meridian $\{x\}\times \partial D^2$ to a circle of slope $p/q \in \Q$
where we use the convention that $ \{x\} \times \partial D^2$ has slope
$\infty$ and $S^1\times \{y\}$ has slope 0.  Thus $L_{1/0} = S^1\times
S^2$ and $L_{0/1} = S^3$.  Recall also that the fraction $p/q$
determines $L_{p/q}$ completely.  At a loss of some precision
Theorem~\ref{main} can be generalized as follows for arbitrary
lens spaces.

\begin{theorem}\label{main_lense}
There exists on a lens space $L_{p/q}$ an overtwisted contact
structure $\zeta$ with the property that the group of contactomorphisms
of $(L_{p/q}, \zeta)$ contains countably many non-conjugate 2-tori.
\end{theorem}

\begin{proof}
The proof is essentially the same as the proof of Theorem~\ref{main}.
We start by giving a  description of lens spaces which is
convenient for our purposes.  Consider again $M= S^1 \times S^1 \times
[0,1]$.  Consider the action of $S^1$ on the boundary of $M$ defined
on $S^1\times S^1 \times \{0\}$ by the vector field
$\frac{\partial}{\partial \theta_2}$ and on $S^1\times S^1 \times
\{1\}$ by the vector field $l\frac{\partial}{\partial \theta_1} -
k\frac{\partial}{\partial \theta_2}$.  By Proposition~\ref{prop2} the cut
$$
X_{k,l} := M/\sim
$$
is a manifold.
Note that the standard action of $ S^1\times S^1$ on $M$
descends to an action on $X_{k,l}$. One can show that $X_{k,l} =
L_{(-k)/l}$.

Next fix $(0,0)\not = (k,l) \in \Z \times \N$ and consider for each
$j\in\N$ the embedding 
$$
\iota_j : M = S^1 \times S^1 \times [0,1] \to S^1 \times S^1 \times \R, 
\quad \iota _j (\theta_1, \theta _2, t) =  
(\theta_1, \theta _2, (\theta_{k,l} + 2\pi j) t)
$$
where $\theta_{k,l}$ is the unique angle with $0<\theta_{k,l} \leq
\pi$ and $\tan \theta_{k,l} = l/k$.  Let $\alpha _j = \iota_j ^*( \cos
t\, d\theta_1 + \sin t \, d\theta_2)$.  By Proposition~\ref{cont_prop2}
each of the contact forms $\alpha_j$ induces a contact from
$\bar{\alpha}_j$ on $X_{k, l}$.  As in the proof of Theorem~\ref{main}
one shows that
\begin{enumerate}
\item  For all $j,j'\geq 0$ the contact structures $\ker \bar{\alpha}_j$ and 
$\ker \bar{\alpha}_{j'}$ are homotopic as 2-plane fields.

\item  The contact  structures $\ker \bar{\alpha}_j$ are overtwisted for 
all $j>0$, hence are all equivalent. 
\end{enumerate} 
Let us denote the contact structure defined by $ \bar{\alpha}_j$,
$j>0$ by $\zeta$.  This is the contact structure in the statement of
the theorem.

Finally by examining the number of connected components of the fibers
the appropriate moment maps one sees that for any two distinct
integers $j, j' \geq 1$, there is no contactomorphism $\psi: (X_{k,l},
\bar{\alpha}_j) \to (X_{k,l}, \bar{\alpha}_{j'})$ which is $S^1\times
S^1$ equivariant.
\end{proof}

\subsection*{Acknowledgments}  
I thank Emmanuel Giroux for a number of useful conversations, in
particular for his help in proving theorem~\ref{main}.  I thank Ilya
Ustilovsky for helpful suggestions and for reading the draft of this
paper.  I thank Peter Ozsvath for a number of useful conversations.  I
thank the referee for a helpful and very constructive critique of the
paper.

\end{document}